\documentclass[a4paper,12pt,leqno]{article}

\usepackage[vmargin=3.5cm]{geometry}

\usepackage{amssymb,amsfonts,amsmath,amsthm}
\usepackage[utf8]{inputenc}
\usepackage{tikz,relsize} 
\usepackage[cmtip,arrow]{xy}
\usepackage{pb-diagram,pb-xy}

\usepackage[english]{babel}

\usepackage{tikz-cd}
\usetikzlibrary{calc}

\usepackage{hyperref}
\usepackage[noblocks]{authblk}
\usepackage{enumerate}

\usepackage{newfloat}
\DeclareFloatingEnvironment[fileext=lod,name=Diagram,placement=p]{fdiagram}

\newcommand\st{\,|\,}

\def\ns#1{\ensuremath\mathbb{#1}}

\def\Z{\ns{Z}}
\def\Q{\ns{Q}}
\def\R{\ns{R}}
\def\C{\ns{C}}

\usepackage[mathscr]{eucal}

\DeclareMathOperator{\Isom}{Isom}
\DeclareMathOperator{\Aff}{Aff}
\DeclareMathOperator{\Aut}{Aut}
\DeclareMathOperator{\Inn}{Inn}
\DeclareMathOperator{\Out}{Out}
\DeclareMathOperator{\GL}{GL}

\DeclareMathOperator{\Iso}{Iso}

\numberwithin{equation}{section}
\newtheorem{thm}{Theorem}[section]
\newtheorem{lm}[thm]{Lemma}
\newtheorem{cor}[thm]{Corollary}

\theoremstyle{definition}
\newtheorem{df}{Definition}[section]
\newtheorem{ex}{Example}[section]

\theoremstyle{remark}

\usepackage[hang,flushmargin]{footmisc}

\author{M. Ha{\l}enda}
\author{R. Lutowski\thanks{Corresponding author: \url{rafal.lutowski@ug.edu.pl}}}
\affil{Institute of Mathematics, Faculty of Mathematics, Physics and Informatics, University of Gda\'nsk, 
	80-308 Gda\'nsk, Poland.}

\title{Symmetries of complex flat manifolds}

\date{}

\divide\dgARROWLENGTH by 2

\begin{document}
	
\renewcommand{\thefootnote}{\fnsymbol{footnote}} 
\renewcommand{\thefootnote}{\arabic{footnote}}

\maketitle

\begin{abstract}
	In this article we show how to calculate the group of automorphisms of flat K\"ahler manifolds. Moreover we are interested in the problem of classification of such manifolds up to biholomorphism. We consider these problems from two points of view. The first one treats the automorphism group as a subgroup of the group of affine transformations,
	while in the second one we analyze it using automorphisms of complex tori. This leads us to the analogues of the Bieberbach theorems in the complex case. We end with some examples, which in particular show that in general the finiteness of the automorphism group depends not only on the fundamental group of a flat manifold.
\end{abstract}

\section{Introduction}

Let $\Gamma$ be a \emph{Bieberbach group} of dimension $n$, i.e. a discrete, cocompact and torsion-free subgroup of the group $\Isom(\R^n) = O(n) \ltimes \R^n$ of isometries of the Euclidean space $\R^n$. Let $L \subset \R^n$ be a group of vectors which define all possible translations in $\Gamma$. We have the following short exact sequence
\begin{equation}
	\label{eq:ses}
	0 \longrightarrow L \stackrel{i}{\longrightarrow} \Gamma \stackrel{r}{\longrightarrow} G \longrightarrow 1,
\end{equation}
where
\[
i(l) = (1,l) \text{ and } r(A,a) = A
\]
for every $l \in L, (A,a) \in \Gamma$. Recall that by the first Bieberbach theorem, $L$ is a free abelian group of rank $n$, it spans the space $\R^n$ and $G$ is a finite group. Moreover the $G$-module structure 
on $L$ is defined by the \emph{holonomy representation of $\Gamma$} as follows:
\[
\varphi(g)(l) = i^{-1}(\gamma i(l)\gamma^{-1}),
\]
where $g \in G, l \in L$ and $r(\gamma) = g$. Let the extension (\ref{eq:ses}) correspond to the cohomology class $\alpha\in H^2(G,L)$. Such a class (i.e. corresponding to a torsion-free extension) is called \emph{special}. In this case the quotient space $\R^n/\Gamma$ is a compact Riemannian manifold $M^n$ of dimension $n$ with
the Riemann curvature tensor equal to zero, a \emph{flat manifold} for short.

In this article we shall consider flat manifolds with a complex structure, i.e. \emph{flat K\"ahler manifolds}.
Recall that $\R^n/\Gamma$ admits a complex structure if and only if $\varphi$ is essentially complex, i.e. $n$ is an even number and each $\R$-irreducible
summand of $\varphi$ which is also $\C$-irreducible occurs with even multiplicity, 
see \cite[Proposition 7.2]{S}. In this case we can consider the group $\Gamma$ as a discrete subgroup of 
$U(\frac{n}{2})\ltimes\C^{\frac{n}{2}}$ (see \cite[Proposition 7.1 (iii)]{S}). 
Equivalently, $L$ admits a complex structure, i.e. a matrix $J \in \GL(n,\R)$ s.t. $J^2=-1$ and $J \varphi(g)=\varphi(g) J$ for every $g \in G$ (see \cite[Proposition 3.1]{J90}).

Let $M\simeq\R^{2n}/\Gamma\simeq\C^{n}/\Gamma$ be a flat K\"ahler manifold of a complex dimension $n$
with a fundamental group $\Gamma$. 
From \cite[Chapter 5]{S} we know that $\Aff(M)/\Aff_{0}(M) \cong \Out(\Gamma)$,
where $\Aff(M)$ is the group of affine self-equivalences of the manifold $M$ with the identity component $\Aff_{0}(M)$ being a torus of the dimension equal to the first Betti number of $M$.
By \cite[page 66]{S} we can define the group $\Aff(M)$ as follows:
\begin{equation}
	\Aff(M) = \{f:M\to M\mid \tilde{f}\in \GL(2n,\R)\ltimes\R^{2n}\}, 
\end{equation}
where $\tilde{f}$ denotes a lift of $f$ in the universal covering space $\R^{2n}$ of $M$.

In our paper we shall consider automorphisms (i.e. biholomorphic self-maps) of flat K\"ahler manifolds. First, observe that if $(M,J)$ is a flat manifold $M$ with a complex structure $J$, then by definition any biholomorphic map $f:M\to M$ satisfies $df \circ J = J\circ df$. Moreover, \cite[Proposition 2.1]{BH} implies that $f$ lifts to an affine map $\tilde{f}$. This justifies the following definition:
\begin{df}
	\label{df:automorphisms}
	Let $M$ be a flat K\"ahler manifold with a complex structure $J$. Then 
	\[
	\Aut(M) = 
	\{f\in \Aff(M)\mid d{f} \circ J = J\circ d{f}\}.
	\]
\end{df}

The paper consists essentially of two parts. In the first part we consider the case when we start with a flat manifold and equip it with a complex structure. Section \ref{sec:basic_diagram} deals with complex analogues of two so-called 9-diagrams for affine self-equivalences of flat manifolds and outer automorphisms of their fundamental groups, described in \cite{S} and \cite{Ch86}. In Section \ref{sec:classification_problem} we extend well known classification criteria of flat manifolds to the case of K\"ahler flat manifolds.

The second part of the paper begins with Section \ref{sec:tori}, in which  results from the previous sections are reformulated. Any flat K\"ahler manifold is a quotient of a complex torus by a finite group acting freely. From this point of view, we are able to describe automorphism group as well as to find biholomorphism classes in terms of some subgroups of the automorphism group of a complex torus. In the conclusion we formulate complex analogues of Bieberbach theorems. We make use of derived results in the next section. It is devoted to some examples, showing differences between groups of automorphisms and affine transformations. 

\section{Basic diagram}

\label{sec:basic_diagram}

Let $\Gamma$ be a Bieberbach group which fits into the short exact sequence \eqref{eq:ses}.
Let $\GL(L)$ be the group of those real matrices which define automorphisms of $L$:
\[
\GL(L) := \{ A \in \GL(n,\R) \st AL = A^{-1}L = L \}.
\]
In other words, $\GL(L)=B \GL(n,\Z) B^{-1}$, where $B \in \GL(n,\R)$ is a matrix whose columns form a basis of $L$. Let $N(G)$ be the normalizer of $G$ in the group $\GL(L)$:
\[
N(G) := N_{\GL(L)}(G) = \{ A \in \GL(L) \st AGA^{-1} = G \}.
\]
Note that we have identified $G$ with its image $\varphi(G) \subset \GL(L)$. 
Let $N_\alpha$ denote the stabilizer of $\alpha \in H^2(G,L)$ -- which corresponds to the extension \eqref{eq:ses} -- under the action $*$ which is induced from the action on cocycles as follows:
\begin{equation}
	\label{eq:star_action}
	n*f(a,b) := nf(n^{-1}an,n^{-1}bn),
\end{equation}
where $n \in N(G), a,b \in G$ and $f\colon G \times G \to L$ is a $2$-cocycle (see
\cite[page 65]{S}). Then by \cite[Theorem 5.1]{S} we have commutative Diagram \ref{diag:basic} with exact rows and columns, where $\Aut^0(\Gamma)$ is a group of all automorphisms of $\Gamma$ which can be defined by the conjugation by some pure translation of $\R^n$.

\begin{fdiagram}
	\[
	\begin{diagram}
		\node{}           \node{0}\arrow{s} \node{1}\arrow{s} \node{1}\arrow{s} \\
		\node{0}\arrow{e} \node{L/L^G}\arrow{s}\arrow{e}\node{\Inn(\Gamma)}\arrow{s}\arrow{e}\node{G}\arrow{s}\arrow{e}\node{1}\\
		\node{1}\arrow{e}\node{\Aut^0(\Gamma)} \arrow{s}\arrow{e} \node{\Aut(\Gamma)}\arrow{s}\arrow{e}\node{N_\alpha}\arrow{s}\arrow{e}\node{1}\\
		\node{0}\arrow{e}\node{H^1(G,L)}\arrow{s}\arrow{e}\node{\Out(\Gamma)}\arrow{s}\arrow{e}\node{N_\alpha/G}\arrow{s}\arrow{e}\node{1}\\
		\node{}\node{0}\node{1}\node{1}
	\end{diagram}
	\]
	\caption{Basic diagram for Bieberbach groups}
	\label{diag:basic}
\end{fdiagram}

Let $M=\R^n/\Gamma$ be the flat manifold with the fundamental group $\Gamma$. By the results presented in \cite[Section V.6]{Ch86} we have commutative Diagram \ref{diag:aff} with exact rows and columns,
\begin{fdiagram}
	\[
	\begin{diagram}
		\node{}           \node{0}\arrow{s} \node{1}\arrow{s} \node{1}\arrow{s} \\
		\node{0}\arrow{e} \node{L^G}\arrow{s}\arrow{e}\node{\Gamma}\arrow{s}\arrow{e}\node{\Inn(\Gamma)}\arrow{s}\arrow{e}\node{1}\\
		\node{0}\arrow{e}\node{(\R^n)^G} \arrow{s}\arrow{e,t}{i}\node{N(\Gamma)}\arrow{s}\arrow{e,t}{\Phi}\node{\Aut(\Gamma)}\arrow{s}\arrow{e}\node{1}\\
		\node{1}\arrow{e}\node{\Aff_0(M)}\arrow{s}\arrow{e}\node{\Aff(M)}\arrow{s}\arrow{e}\node{\Out(\Gamma)}\arrow{s}\arrow{e}\node{1}\\
		\node{}\node{1}\node{1}\node{1}
	\end{diagram}
	\]
	\caption{Affine self-equivalences of flat manifolds}
	\label{diag:aff}
\end{fdiagram}
where $N(\Gamma)$ denotes the normalizer of $\Gamma$ in the group $\Aff(\R^n) = \GL(n,\R) \ltimes \R^n$, the groups $L^G$ and $(\R^n)^G$ are the groups of fixed points under the action (left multiplication) of $G$ on $L$ and $\R^n$ respectively. Let us note some basic facts:
\begin{lm}
	Let $(A,a),(B,b) \in N(\Gamma)$. Then:
	\begin{enumerate}
		\item If $(A,a)$ and $(B,b)$ define the same automorphism of $\Gamma$ then $A = B$.
		\item If $(A,a)$ and $(B,b)$ are lifts of the same affine self-equivalence of $M$ then $A^{-1}B \in G$.
		\item $N_\alpha = r(N(\Gamma)) = \{ d \tilde f \st \tilde f \in N(\Gamma) \}$, i.e. it is a group of differentials of all possible lifts of elements of
		$\Aff(M)$.
	\end{enumerate}
\end{lm}

The proof of the above lemma takes advantage of diagrams \ref{diag:basic} and \ref{diag:aff} as well as the description of differential structures of flat manifolds presented in \cite[pp. 50-52]{Ch86}.

Now assume that $n$ is even and that a matrix $J \in \GL(n,\R)$ defines a complex structure on $M$. Let $f \in \Aff(M)$ and $(A,a)$ be a lift of $f$. We get that
\[
f \in \Aut(M) \Leftrightarrow JA = AJ.
\]

Define
\[
N^J(\Gamma) := \{  (A,a) \in N(\Gamma) \st AJ = JA \}.
\]
We immediately get
\begin{cor} \ 
	\begin{enumerate}
		\item $\Gamma \subset N^J(\Gamma)$.
		\item $i( (\R^n)^G ) \subset N^J(\Gamma)$.
		\item $\Aut^0(\Gamma) \subset \Phi( N^J(\Gamma) )$.
		\item $N^J_\alpha := r(N^J(\Gamma)) = \{ A \in N_\alpha \st AJ = JA \}$.
	\end{enumerate}
\end{cor}

Denote $\Aut^J(\Gamma) = \Phi( N^J(\Gamma) ), \Out^J(\Gamma) = \Aut^J(\Gamma) / \Inn(\Gamma)$. Diagrams \ref{diag:basic_c} and \ref{diag:bihol} are the analogs of diagrams \ref{diag:basic} and \ref{diag:aff} in the complex case.

At the end of the section let us note that we can link diagrams \ref{diag:basic} with \ref{diag:aff} and \ref{diag:basic_c} with \ref{diag:bihol} and get commutative diagrams \ref{diag:cube_real} and \ref{diag:cube_complex} respectively, by doing quite easy diagram chase. Helpful, for understanding certain part of the diagrams, description of $H^1(G,L)$ may be found in \cite[Section 4.2]{Lu13}. Note that for every line in both diagrams the first arrow is a monomorphism and the second one -- epimorphism. The groups which were not yet defined include:
\begin{itemize}
	\item $Z(\Gamma)$ -- the center of $\Gamma$,
	\item $N^0(\Gamma)$ -- the group of all possible translations in $N(\Gamma)$,
	\item $C(\Gamma) \cong (\R^n)^G$ -- the centralizer of $\Gamma$ in $\Aff(\R^n)$,
	\item $T^G$ -- the group of fixed points of the action of $G$ on the torus $T=\R^n/\Z^n$.
\end{itemize}

\begin{fdiagram}
	\[
	\begin{diagram}
		\node{}           \node{0}\arrow{s} \node{1}\arrow{s} \node{1}\arrow{s} \\
		\node{0}\arrow{e} \node{L/L^G}\arrow{s}\arrow{e}\node{\Inn(\Gamma)}\arrow{s}\arrow{e}\node{G}\arrow{s}\arrow{e}\node{1}\\
		\node{1}\arrow{e}\node{\Aut^0(\Gamma)} \arrow{s}\arrow{e} \node{\Aut^J(\Gamma)}\arrow{s}\arrow{e}\node{N^J_\alpha}\arrow{s}\arrow{e}\node{1}\\
		\node{0}\arrow{e}\node{H^1(G,L)}\arrow{s}\arrow{e}\node{\Out^J(\Gamma)}\arrow{s}\arrow{e}\node{N^J_\alpha/G}\arrow{s}\arrow{e}\node{1}\\
		\node{}\node{0}\node{1}\node{1}
	\end{diagram}
	\]
	\caption{Basic diagram for Bieberbach group with an almost complex structure}
	\label{diag:basic_c}
\end{fdiagram}

\begin{fdiagram}
	\[
	\begin{diagram}
		\node{}           \node{0}\arrow{s} \node{1}\arrow{s} \node{1}\arrow{s} \\
		\node{0}\arrow{e} \node{L^G}\arrow{s}\arrow{e}\node{\Gamma}\arrow{s}\arrow{e}\node{\Inn(\Gamma)}\arrow{s}\arrow{e}\node{1}\\
		\node{0}\arrow{e}\node{(\R^n)^G} \arrow{s}\arrow{e}\node{N^J(\Gamma)}\arrow{s}\arrow{e}\node{\Aut^J(\Gamma)}\arrow{s}\arrow{e}\node{1}\\
		\node{1}\arrow{e}\node{\Aff_0(M)}\arrow{s}\arrow{e}\node{\Aut(M)}\arrow{s}\arrow{e}\node{\Out^J(\Gamma)}\arrow{s}\arrow{e}\node{1}\\
		\node{}\node{1}\node{1}\node{1}
	\end{diagram}
	\]
	\caption{Automorphisms of complex flat manifolds}
	\label{diag:bihol}
\end{fdiagram}

\begin{fdiagram}
	\begin{tikzpicture}
		\matrix[matrix of nodes, name=m, every node/.style={minimum width=3em}, commutative diagrams/row sep=1em, ampersand replacement=\&] {
			\& \&$L^G$ \& \& \& $L$  \& \& \& $L/L^G$\\
			\&$Z(\Gamma)$\& \& \& $\Gamma$ \& \& \& $\Inn(\Gamma)$ \\
			$1$\&\&\& $G$ \&\&\& $G$\\
			\& \&$(\R^n)^G$  \& \& \& $N^0(\Gamma)$  \& \& \& $\Aut^0(\Gamma)$ \\
			\&$C(\Gamma)$\& \& \& $N(\Gamma)$ \& \& \& $\Aut(\Gamma)$\\
			$1$\& \& \& $N_\alpha$ \& \& \& $N_\alpha$\\
			\& \& $\Aff_0(X)$ \& \& \& $T^G$  \& \& \& $H^1(G,L)$\\
			\&$\Aff_0(M)$\& \& \& $\Aff(M)$ \& \& \& $\Out(\Gamma)$ \\
			$1$\& \& \& $N_\alpha/G$  \& \& \& $N_\alpha/G$\\
		};
		\path[commutative diagrams/.cd, every arrow, every label]
		(m-1-3) edge (m-1-6) 
		(m-1-3) edge (m-4-3)
		(m-1-3) edge (m-2-2)
		(m-1-6) edge (m-1-9)
		(m-1-6) edge (m-4-6)
		(m-1-6) edge (m-2-5)
		(m-1-9) edge (m-4-9)
		(m-1-9) edge (m-2-8)
		(m-4-3) edge (m-4-6)
		(m-4-3) edge (m-7-3)
		(m-4-3) edge (m-5-2)
		(m-4-6) edge (m-4-9)
		(m-4-6) edge (m-7-6)
		(m-4-6) edge (m-5-5)
		(m-4-9) edge (m-7-9)
		(m-4-9) edge (m-5-8)
		(m-7-3) edge (m-7-6)
		(m-7-3) edge (m-8-2)
		(m-7-6) edge (m-7-9)
		(m-7-6) edge (m-8-5)
		(m-7-9) edge (m-8-8)
		
		(m-2-2) edge[commutative diagrams/crossing over] (m-2-5)
		(m-2-2) edge[commutative diagrams/crossing over] (m-5-2)
		(m-2-2) edge[commutative diagrams/crossing over] (m-3-1)
		(m-2-5) edge[commutative diagrams/crossing over] (m-2-8)
		(m-2-5) edge[commutative diagrams/crossing over] (m-5-5)
		(m-2-5) edge[commutative diagrams/crossing over] (m-3-4)
		(m-2-8) edge[commutative diagrams/crossing over] (m-5-8)
		(m-2-8) edge[commutative diagrams/crossing over] (m-3-7)
		(m-5-2) edge[commutative diagrams/crossing over] (m-5-5)
		(m-5-2) edge[commutative diagrams/crossing over] (m-8-2)
		(m-5-2) edge[commutative diagrams/crossing over] (m-6-1)
		(m-5-5) edge[commutative diagrams/crossing over] (m-5-8)
		(m-5-5) edge[commutative diagrams/crossing over] (m-8-5)
		(m-5-5) edge[commutative diagrams/crossing over] (m-6-4)
		(m-5-8) edge[commutative diagrams/crossing over] (m-8-8)
		(m-5-8) edge[commutative diagrams/crossing over] (m-6-7)
		(m-8-2) edge[commutative diagrams/crossing over] (m-9-1)
		(m-8-2) edge[commutative diagrams/crossing over] (m-8-5)
		(m-8-5) edge[commutative diagrams/crossing over] (m-9-4)
		(m-8-5) edge[commutative diagrams/crossing over] (m-8-8)
		(m-8-8) edge[commutative diagrams/crossing over] (m-9-7)
		(m-3-1) edge[commutative diagrams/crossing over] (m-3-4)
		(m-3-1) edge[commutative diagrams/crossing over] (m-6-1)
		(m-3-4) edge[commutative diagrams/crossing over] (m-3-7)
		(m-3-4) edge[commutative diagrams/crossing over] (m-6-4)
		(m-3-7) edge[commutative diagrams/crossing over] (m-6-7)
		(m-6-1) edge[commutative diagrams/crossing over] (m-6-4)
		(m-6-1) edge[commutative diagrams/crossing over] (m-9-1)
		(m-6-4) edge[commutative diagrams/crossing over] (m-6-7)
		(m-6-4) edge[commutative diagrams/crossing over] (m-9-4)
		(m-6-7) edge[commutative diagrams/crossing over] (m-9-7)
		(m-9-1) edge[commutative diagrams/crossing over] (m-9-4)
		(m-9-4) edge[commutative diagrams/crossing over] (m-9-7)
		;
	\end{tikzpicture}
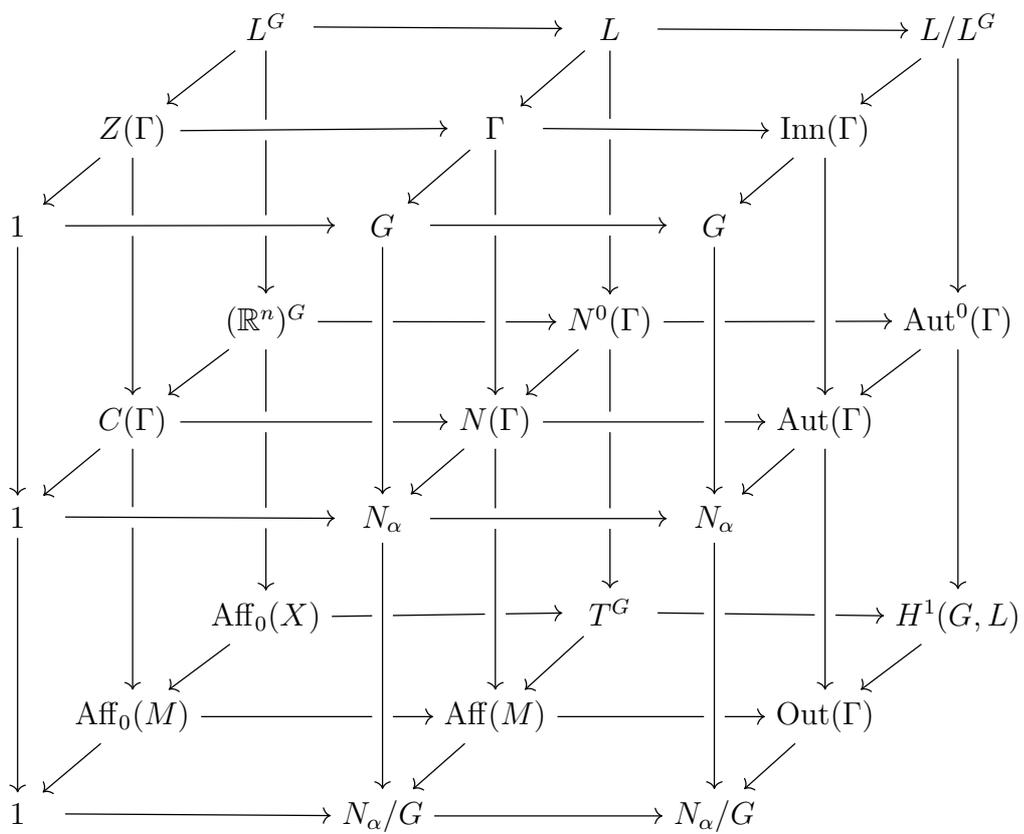
	\caption{,,Cube'' diagram in the real case}
	\label{diag:cube_real}
\end{fdiagram}

\begin{fdiagram}
	\begin{tikzpicture}
		\matrix[matrix of nodes, name=m, every node/.style={minimum width=3em}, commutative diagrams/row sep=1em, ampersand replacement=\&] {
			\& \&$L^G$ \& \& \& $L$  \& \& \& $L/L^G$\\
			\&$Z(\Gamma)$\& \& \& $\Gamma$ \& \& \& $\Inn(\Gamma)$ \\
			$1$\&\&\& $G$ \&\&\& $G$\\
			\& \&$(\R^n)^G$  \& \& \& $N^0(\Gamma)$  \& \& \& $\Aut^0(\Gamma)$ \\
			\&$C(\Gamma)$\& \& \& $N^J(\Gamma)$ \& \& \& $\Aut^J(\Gamma)$\\
			$1$\& \& \& $N^J_\alpha$ \& \& \& $N^J_\alpha$\\
			\& \& $\Aff_0(X)$ \& \& \& $T^G$  \& \& \& $H^1(G,L)$\\
			\&$\Aff_0(M)$\& \& \& $\Aut(M)$ \& \& \& $\Out^J(\Gamma)$ \\
			$1$\& \& \& $N^J_\alpha/G$  \& \& \& $N^J_\alpha/G$\\
		};
		\path[commutative diagrams/.cd, every arrow, every label]
		(m-1-3) edge (m-1-6) 
		(m-1-3) edge (m-4-3)
		(m-1-3) edge (m-2-2)
		(m-1-6) edge (m-1-9)
		(m-1-6) edge (m-4-6)
		(m-1-6) edge (m-2-5)
		(m-1-9) edge (m-4-9)
		(m-1-9) edge (m-2-8)
		(m-4-3) edge (m-4-6)
		(m-4-3) edge (m-7-3)
		(m-4-3) edge (m-5-2)
		(m-4-6) edge (m-4-9)
		(m-4-6) edge (m-7-6)
		(m-4-6) edge (m-5-5)
		(m-4-9) edge (m-7-9)
		(m-4-9) edge (m-5-8)
		(m-7-3) edge (m-7-6)
		(m-7-3) edge (m-8-2)
		(m-7-6) edge (m-7-9)
		(m-7-6) edge (m-8-5)
		(m-7-9) edge (m-8-8)
		
		(m-2-2) edge[commutative diagrams/crossing over] (m-2-5)
		(m-2-2) edge[commutative diagrams/crossing over] (m-5-2)
		(m-2-2) edge[commutative diagrams/crossing over] (m-3-1)
		(m-2-5) edge[commutative diagrams/crossing over] (m-2-8)
		(m-2-5) edge[commutative diagrams/crossing over] (m-5-5)
		(m-2-5) edge[commutative diagrams/crossing over] (m-3-4)
		(m-2-8) edge[commutative diagrams/crossing over] (m-5-8)
		(m-2-8) edge[commutative diagrams/crossing over] (m-3-7)
		(m-5-2) edge[commutative diagrams/crossing over] (m-5-5)
		(m-5-2) edge[commutative diagrams/crossing over] (m-8-2)
		(m-5-2) edge[commutative diagrams/crossing over] (m-6-1)
		(m-5-5) edge[commutative diagrams/crossing over] (m-5-8)
		(m-5-5) edge[commutative diagrams/crossing over] (m-8-5)
		(m-5-5) edge[commutative diagrams/crossing over] (m-6-4)
		(m-5-8) edge[commutative diagrams/crossing over] (m-8-8)
		(m-5-8) edge[commutative diagrams/crossing over] (m-6-7)
		(m-8-2) edge[commutative diagrams/crossing over] (m-9-1)
		(m-8-2) edge[commutative diagrams/crossing over] (m-8-5)
		(m-8-5) edge[commutative diagrams/crossing over] (m-9-4)
		(m-8-5) edge[commutative diagrams/crossing over] (m-8-8)
		(m-8-8) edge[commutative diagrams/crossing over] (m-9-7)
		(m-3-1) edge[commutative diagrams/crossing over] (m-3-4)
		(m-3-1) edge[commutative diagrams/crossing over] (m-6-1)
		(m-3-4) edge[commutative diagrams/crossing over] (m-3-7)
		(m-3-4) edge[commutative diagrams/crossing over] (m-6-4)
		(m-3-7) edge[commutative diagrams/crossing over] (m-6-7)
		(m-6-1) edge[commutative diagrams/crossing over] (m-6-4)
		(m-6-1) edge[commutative diagrams/crossing over] (m-9-1)
		(m-6-4) edge[commutative diagrams/crossing over] (m-6-7)
		(m-6-4) edge[commutative diagrams/crossing over] (m-9-4)
		(m-6-7) edge[commutative diagrams/crossing over] (m-9-7)
		(m-9-1) edge[commutative diagrams/crossing over] (m-9-4)
		(m-9-4) edge[commutative diagrams/crossing over] (m-9-7)
		;
	\end{tikzpicture}
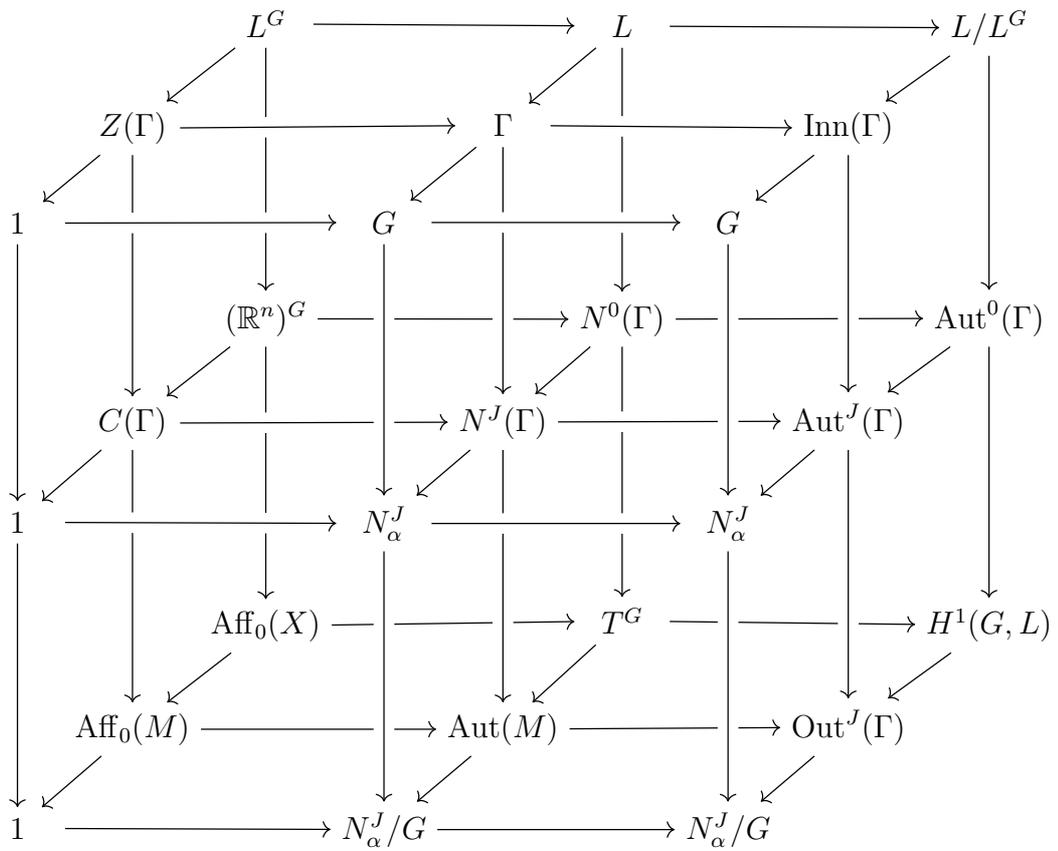
	\caption{,,Cube'' diagram in the complex case}
	\label{diag:cube_complex}
\end{fdiagram}

\section{Classification problem}

\label{sec:classification_problem}

Let the pairs $(\Gamma,J)$ and $(\Gamma',J')$ define flat K\"ahler manifolds, i.e. $\Gamma,\Gamma' \subset \Iso(\R^n)$ are Bieberbach groups and $J,J' \in \GL(n,\R)$ define complex structures on the flat manifolds $M=\R^n/\Gamma, M'=\R^n/\Gamma'$ respectively, where $n$ is an even integer. 

Assume that the manifolds $M$ and $M'$ are biholomorphic and let $f \colon M \to M'$ be any biholomorphism. The groups $\Gamma$ and $\Gamma'$ are isomorphic and if $\tilde{f} = (A,a) \in \Aff(\R^n)$ is a lift of $f$ then the isomorphism $\phi \colon \Gamma \to \Gamma'$ is given by the formula 
\[
\gamma \mapsto (A,a) \gamma (A,a)^{-1}.
\]
Moreover in this case we get
\[
A J = J' A.
\]

The above considerations show that in the classification problem of biholomorphic flat K\"ahler manifolds it is enough to look at the extensions of the form \eqref{eq:ses} with the same complex structure $J$. Let $\Gamma$ and $\Gamma'$ be extensions of $L$ by $G$ as in \eqref{eq:ses}, defined by the cohomology classes $\alpha, \alpha' \in H^2(G,L)$ respectively. Recall that $N(G)$, the normalizer of $G$ in $\GL(L)$, acts by $*$ on $H^2(G,L)$ (see \eqref{eq:star_action}). We have the following theorem:

\begin{thm}
	The flat K\"ahler manifolds $M=\R^n/\Gamma$ and $M'=\R^n/\Gamma'$, both with the same almost complex structure $J$, are biholomorphic if and only if there exists $A \in N(G)$ such that
	\[
	A*\alpha = \alpha' \text{ and } AJ = JA.
	\]
\end{thm}

\begin{proof}
	Let $f \colon M \to M'$ be a biholomorphism and let $\tilde{f} = (A,a) \in \Aff(\R^n)$ be its lift. By the above description of the isomorphism of $\Gamma$ and $\Gamma'$ we have that $A * \alpha = \alpha'$. Now, since $A = df$, hence $AJ = JA$.
	
	On the other hand if there exists $A \in N(G)$ st. $AJ=JA$ and $A * \alpha = \alpha'$ then we can find $a \in \R^n$ st. conjugation by $(A,a)$ defines an isomorphism of $\Gamma$ and $\Gamma'$. In this case $(A,a)$ is a lift (in $\R^n$) of an affine map $f \colon M \to M'$, which by Definition \ref{df:automorphisms} must be a biholomorphism.
\end{proof}

\begin{cor}
	\label{cor:bihol_classes}
	The set of biholomorphism classes of flat K\"ahler manifolds of the form $\R^n/\Gamma$ with the almost complex structure $J$ and for which the group $\Gamma$ fits into the short exact sequence \eqref{eq:ses} is in one to one correspondence with the orbits of the action $*$ of the group
	\[
	N^J = \{ A \in N(G) \st AJ = JA \}
	\]
	on the set of special cohomology classes of $H^2(G,L)$.
\end{cor}

\section{Automorphisms of quotients of complex tori}
\label{sec:tori}
By the results of previous sections it is possible to compute biholomorphism group of a flat K\"ahler manifold (with fixed complex structure $J$) or to find all biholomorphism classes of flat K\"ahler manifolds with the same fundamental group $\Gamma$ and fixed complex structure $J$. In some cases these computations can be done with the help of computer package CARAT (\cite{carat}). However, we will rephrase our results in the following way.

Assume, that a flat K\"ahler manifold $M$ is defined -- as in the previous section -- by a pair $(\Gamma,J)$ where $\Gamma \subset \Iso(\R^{2n})$.
Since $\Gamma$ fits into short exact sequence 
\[
0\longrightarrow L \longrightarrow \Gamma\longrightarrow G\longrightarrow1,
\]
then $T=\R^{2n}/L$ is a complex torus (with a complex structure $J$) and $M=T/\tilde{G}$, where $\tilde{G}$ is a subgroup of the group $\Aut(T)$, isomorphic to $G$. Recall that if $T$ is a complex torus, then we have the following short exact sequence:
\[
0\longrightarrow T\longrightarrow \Aut(T)\stackrel{\pi}{\longrightarrow} \Aut_0(T)\longrightarrow 1
\]
where $\Aut_0(T)$ is the group of all biholomorphic homomorphisms of $T$ (see \cite[Proposition 2.1]{BH}. Since $L$ is maximal 
then group $\tilde{G}$ may be identified with a class $\alpha\in H^1(G,T)\cong H^2(G,L)$, where we regard $G$ as a subgroup $\pi(\tilde G)$ of $\Aut_0(T)$. 

\begin{df}
	Let $M=T/\tilde{G}$ be a flat K\"ahler manifold and assume that $\tilde{G}$ corresponds to a cohomology class $\alpha\in H^1(G,T)$. Then we define the following groups:
	\begin{enumerate}
		\item $N_{\Aut_0(T)}(G)$ -- the normalizer of $G$ in $\Aut_0(T)$,
		\item $N_\alpha$ -- the stabilizer of $\alpha \in H^1(G,T)$ under the action $*$ of $N_{\Aut_0(T)}(G)$, defined similarly as in \eqref{eq:star_action}, that is:
		\[
		n*f(a) := n f(n^{-1}an),
		\]
		where $a \in G, n \in N_{\Aut_0(T)}(G)$ and $f \colon G \to T$ is one-cocycle.
	\end{enumerate}
\end{df}
Note that $N_\alpha$ defined above differ from the one defined in Section \ref{sec:basic_diagram} and is in fact isomorphic to the group $N_\alpha^J$ (for an appropriate complex structure $J$).

Corollary \ref{cor:bihol_classes} may be then rephrased in the following way:
\begin{thm}
	Let $T$ be a complex torus and $G\subset \Aut_0(T)$ be a finite group. Then 
	biholomorphism classes of flat K\"ahler manifolds of the form $T/\tilde{G}$ correspond bijectively to the orbits of special cohomology classes of $H^1(G,T)$ under the action $*$ of $N_{\Aut_0(T)}(G)$.
\end{thm}

Moreover, an immediate consequence of Diagram \ref{diag:cube_complex} is:

\begin{thm}
	Let $T$ be a complex torus, $G\subset \Aut_0(T)$ be a finite group and $\alpha\in H^1(G,T)$ a special cohomology class. If $M$ is a flat K\"ahler manifold corresponding to $\alpha$, then the group $\Aut(M)$ of automorphisms of $M$ is an extension of the form
	\[
	0\longrightarrow T^G\longrightarrow \Aut(M)\longrightarrow N_\alpha/G\longrightarrow 1.
	\]
\end{thm}
Now, similarly as in the real case, since $H^1(G,T)$ is finite, by the orbit-stabilizer theorem we get:
\begin{cor}
	Let $M=T/\tilde{G}$ be a flat K\"ahler manifold. The group of automorphisms $\Aut(M)$ is finite if and only if the first Betti number of $M$ is equal to zero
	and the group $N_{\Aut_0(T)}(G)$ is finite.
\end{cor}

Summing up, we may state analogues of Bieberbach theorems for flat K\"ahler manifolds:
\begin{thm}[Complex Bieberbach Theorems]\mbox{}
	\begin{enumerate}
		\item $M$ is a flat K\"ahler manifold of complex dimension $n$ if and only if there exists complex $n$-dimensional torus $T$ and a finite group $\tilde{G}\subset \Aut(T)$ acting freely on $T$ such that $M=T/\tilde{G}$.
		\item Two flat K\"ahler manifolds $M=T/\tilde{G}$ and $M'=T'/\tilde{G}'$ are biholomorphic if and only if there exists a biholomorphic map $\varphi:T\to T'$ such that $\tilde{G} = \varphi^{-1}\tilde{G}'\varphi$.
		\item For every complex torus $T$ there exist only finite number of flat K\"ahler ma\-ni\-folds of the form $T/\tilde{G}$, up to biholomorphism.
	\end{enumerate}
\end{thm}

\section{Examples}

The goal of this section is to look at some examples, which will be used to illustrate two matters. The first one is a classification of complex flat manifolds whose fundamental group is a given Bieberbach group. 
This is somehow a starting point to the second problem, which is to determine automorphisms groups of complex manifolds obtained in the previous step, and to compare these groups with 
the affine self-equivalences group of the underlying real manifold.

First of these examples, presented in subsection \ref{sec:c3xc3}, corresponds to the group that has been constructed by Hiller and Sah in \cite[Proposition 3.3]{HS}. In this case complex manifold structure is unique. Moreover real -- and hence complex -- symmetry groups are both finite. Complex Hantzsche-Wendt manifolds from subsection \ref{sec:chw} (see \cite{H}) are on the other hand examples of manifolds with infinite real, but finite complex symmetries (regardless of a particular complex structure; however exact symmetry type may vary). In the last subsection we discuss a connection between the finiteness of the automorphism group and the holonomy representation, which seems to be much more subtle than in the real case.

\subsection{Complex flat fourfold with holonomy \texorpdfstring{$\Z_3^2$}{Z32}.}
\label{sec:c3xc3}
Let us consider an abelian variety $T=E_\xi\times E_\xi\times E_\xi\times E_\xi$, where $\xi=e^{2\pi i/3}$ and $E_\xi$ denotes the elliptic curve $\mathbb{C}/(\mathbb{Z}+\xi\mathbb{Z})$. Let $G<\Aut_0(T)$ be the group generated by $g_1=(1,\xi,\xi,\xi)$ and $g_2=(\xi,1,\xi,\xi^2)$. We shall show that there exists an element of $H^1(G,T)$ which defines fixed-point free action of a group $\tilde{G} \cong \Z_3^2$ on $T$. 

Let $z\in Z^1(G,T)$ be a 1-cocycle. Adding a suitable coboundary we may assume that $z(g_1)=(a,0,0,0)$ and $z(g_2)=(0,b,c,d)$. Since every element of $G$ has order 3, then the points $a,b,c,d$ are 3-torsion points of $T$. Moreover since $G$ is abelian, then we get that the points $a,b,c,d$ satisfy an equation $\xi x=x$. There are exactly three 3-torsion points with that property and we have to choose non-zero points, since the group of biholomorphisms of $T$ defined by $z$ has to act freely.  Thus $a,b,c,d\in\{\pm\frac{1}{3}(1-\xi)\}$.

Fix a choice of $a,b,c,d$ and let $\Gamma$ be the fundamental group
of $M=T/\tilde{G}$, where $\tilde{G}$ is an extension defined by $z$. The integral holonomy representation of $\Gamma$ is given by matrices:
\[
\varrho(g_1)=
\begin{bmatrix}I&&&\\&\Xi&&\\&&\Xi&\\&&&\Xi\end{bmatrix},\ 
\varrho(g_2)=
\begin{bmatrix}\Xi&&&\\&I&&\\&&\Xi&\\&&&\Xi^2\end{bmatrix},
\]
where $I$ is the $2\times 2$ identity matrix and $\Xi=\begin{bmatrix}0&-1\\1&-1\end{bmatrix}$.  
Observe, that if $M'=T'/\tilde{G}$ is another flat complex manifold with the same integral holonomy, then $T'$ is biholomorphic to $T$. It follows from the fact that this representation determines decomposition of the lattice $L$ which yields decomposition of the torus $T=\mathbb{C}^4/L$ as a product of four elliptic curves, each of them allowing faithful action of $\mathbb{Z}_3$. Moreover the action of the holonomy group is unique up to conjugacy in $\Aut_0(T)$.

Let us compute the normalizer of $G$ in $\Aut_0(T)$. 
Denote the constituents of $T$ by $E_1,E_2,E_3$ and $E_4$
(they are all isomorphic to $E_\xi$). Observe, that if $n\in N_{\Aut_0(T)}(G)$ is a permutation of the curves $E_i$, then it fixes $E_1$ and it may be any even permutation of the remaining curves. This follows from the fact, that $g_1$ cannot be conjugated in $\Aut_0(T)$ to any other element of $G$, and $E_1$ is the only curve on which $g_1$ acts trivially. The evenness comes from the fact that $n g_2 n^{-1} \in G$. On the other hand, if $n\in N_{\Aut_0(T)}(G)$ fixes all curves $E_i$, then it may act on any of them by an element of the group $\langle-\xi\rangle\simeq \mathbb{Z}_6$. Thus $N_{\Aut_0(T)}(G)\cong \mathbb{Z}_6^4\rtimes A_3$.

Now it is easy to observe, that under the action of $N_{\Aut_0(T)}(G)$ on $H^1(G,T)$ all cohomology classes corresponding to fixed-point free actions 
are in the same orbit. Summing up, there is only one complex flat manifold $M$ with the integral holonomy representation $\varrho$ up to biholomorphism.

In the next step we compute the group $\Aut(M)$. Since the action of $G$ on $T$ is diagonal, $T^G\cong \mathbb{Z}_3^4$ (each constituent is generated by $\frac{1}{3}(1-\xi)$). 
For a cohomology class $\alpha \in H^1(G,T)$ corresponding to $M$ we take the one defined by the 1-cocycle $z$, where $a=b=c=d=\frac{1}{3}(1-\xi)$.
For this class we have that $N_{\alpha}$ is exactly $\mathbb{Z}_3^4\rtimes A_3$. Thus $N_{\alpha}/G\cong \mathbb{Z}_3^2\rtimes A_3$ and $\Aut(M)$ is an extension of the form
$$ 0\longrightarrow\Z_3^4\longrightarrow \Aut(M)\longrightarrow \mathbb{Z}_3^2\rtimes A_3\longrightarrow 1.$$
The group of all affine automorphisms of underlying real manifold of $M$ is also a finite group. 
One can check, using for example algorithms presented in \cite{Lu13}, that it is of order $34992$
and the index of $\Aut(M)$ in $\Aff(M)$ is equal to $16$.

\subsection{Complex Hantzsche-Wendt threefolds}
\label{sec:chw}
As our next example, we will examine complex Hantzsche-Wendt manifolds (or shortly -- CHW manifolds) of complex dimension 3. By definition, these are complex flat manifolds with holonomy group $\Z_2^2\subset SU(3)$ (see \cite{H}). There are four possible integral holonomy representations. However, for simplicity we will assume that the integral holonomy is diagonal. All CHW threefolds with this property are diffeomorphic, and their fundamental group has CARAT symbol min.185.1.1.21. Still, there are infinitely many such manifolds up to biholomorphism, which follows from the structure theorem of \cite{H}. In our case this theorem simplifies to the following statement: 
\begin{thm}
	Manifold $M$ is a CHW threefold with diagonal integral holonomy representation if and only if $M$ is a orbit space $T/\tilde G$, where $T=E_1\times E_2\times E_3$ is a product of some elliptic curves and $\tilde G$ is a group generated by mappings $\tilde g_1,\tilde g_2:T\to T$ such that:
	\[
	\begin{array}{ccl}
		\tilde g_1(x_1,x_2,x_3)&=&(x_1+a_1,-x_2+a_2,-x_3+a_3),\\ 
		\tilde g_2(x_1,x_2,x_3)&=&(-x_1+b_1,x_2+b_2,-x_3+b_3),
	\end{array}
	\]
	where $a_1$, $b_2$ and $a_3-b_3$ are non-zero $2$-torsion points of the elliptic curves $E_1,E_2$ and $E_3$ respectively.
\end{thm}

Let $M=T/\tilde{G}$ be as above, $G = \pi(\tilde G), g_1=\pi(\tilde g_1)$ and $g_2=\pi(\tilde g_2)$. If $z\in Z^1(G,T)$ is a cocycle associated to $M$, then adding a suitable coboundary we may assume that $z(g_1)=(a,0,0)$ and $z(g_2)=(0,b,c)$, where $a,b,c$ are non-zero 2-torsion points. In contrary to the previous case, we may choose infinitely many different complex structures on $T$ (using three parameters). Moreover, we shall show that biholomorphism class of $M$ is not uniquely determined by the complex structure on $T$ and a number of non-biholomorphic CHW manifolds covered by the same complex torus $T$ depends on the complex structure of $T$.

Observe, that the group $\Aut_0 (T)$ can be finite (in generic case where $T$ is a product of pairwise non-isomorphic elliptic curves) or infinite. However, the normalizer of $G$ in $\Aut_0 (T)$ is always a finite group. Assume that $n\in N_{\Aut_0(T)}(G)$. Then:
\begin{enumerate}[(i)]
	\item if $n$ centralizes $G$ then it acts on each elliptic curve $E_k$ via an element of $\Aut_0(E_k)$ (a group isomorphic to $\Z_2$, $\Z_4$ or $\Z_6$),
	\item if $n$ does not centralize $G$, then it has to permute elliptic curves $E_1$, $E_2$ and $E_3$ (since complex holonomy representation has three pairwise non-isomorphic constituents). However it is possible only when some of the above elliptic curves are biholomorphic.
\end{enumerate}
Summing up:
\[
N_{\Aut_0(T)}(G)=\left(\Aut_0(E_1)\oplus \Aut_0(E_2)\oplus \Aut_0(E_3)\right)\rtimes S,
\]
where $S$ is group of permutations of a set of those indices $k\in\{1,2,3\}$ for which elliptic curves $E_i$ are biholomorphic.

As we have observed, cohomology classes of the group $H^1(G,T)$ corresponding to torsion-free extensions can be parametrized by $a,b,c$ -- non-zero 2-torsion points of curves $E_1,E_2$ and $E_3$. Let $E_\tau=\C/(\Z+\tau \Z)$ be an elliptic curve. The action of $\Aut_0(E)$ on set of non-zero 2-torsion  points of $E$ has:
\begin{itemize}
	\item one orbit $\left\{\frac{1}{2},\frac{\tau}{2},\frac{1+\tau}{2}\right\}$ if $\tau=e^{\frac{2\pi i}{3}}$,
	\item two orbits $\left\{\frac{1}{2},\frac{\tau}{2}\right\}$ and $\left\{\frac{1+\tau}{2}\right\}$ when $\tau=i$,
	\item three one-element orbits if $E_\tau$ is not biholomorphic to any of the two previous elliptic curves.
\end{itemize}

To calculate the number $m$ of orbits of the action of $N_{\Aut_0 T} (G)$ on special cohomology classes it is now enough to see that if $n\in N_{\Aut_0 T} (G)$ does not centralize $G$, then its action is just an appropriate permutation of points $a$, $b$ and $c$. 
\begin{ex}\label{exCHW}
	Let $T=E_1\times E_2\times E_3$ and $m$ be a number of CHW manifolds with diagonal holonomy covered by $T$.
	\begin{enumerate}
		\item If all elliptic curves $E_1,E_2,E_3$ are biholomorphic to $E_\xi$, where $\xi=e^{\frac{2\pi i}{3}}$, then $m=1$.
		\item If none of elliptic curves $E_1,E_2,E_3$ is biholomorphic to $E_\xi$ or $E_i$ and any two of them are not biholomorphic, then $m=27$.
		\item If $E_1=E_2=E_\xi$ and $E_3=E_i$, then $m=2$. Related two cohomology classes can be parametrized for example by $(a,b,c)=(\frac{1}{2},\frac{1}{2},\frac{1}{2})$ and  $(a',b',c')=(\frac{1}{2},\frac{1}{2},\frac{1+i}{2})$.\label{exCHW3}
	\end{enumerate}
\end{ex}


Finally we will compute $\Aut(M)$. Clearly $T^G=\Z_2^6$ is the set of all $2$-torsion points of $T$. Group $N_\alpha$ depends on particular cohomology class $\alpha\in H^1(G,T)$, but -- as a subgroup of $N_{\Aut_0(T)}(G)$ -- it is finite. Thus $\Aut(M)$ is a finite group. 
\begin{ex}
	Let $M_1$ and $M_2$ be two manifolds from point \ref{exCHW3} of Example~\ref{exCHW}. Then the orders of $\Aut(M_1)$ and $\Aut(M_2)$ are $2^8$ and $2^{9}$ respectively.
\end{ex}
Indeed, we have $N_{\Aut_0(T)}(G)=(\Z_6^2\rtimes S_2)\times \Z_4$. If $\alpha$ corresponds to the values $(a,b,c)=(\frac{1}{2},\frac{1}{2},\frac{1}{2})$ then $N_\alpha=(\Z_2^2\rtimes S_2)\times\Z_2$ (it is generated by transposition on two first elliptic curves and by $-id$ on each of the three curves). Similarly if $\beta$ corresponds to the values $(a',b',c')=(\frac{1}{2},\frac{1}{2},\frac{1+i}{2})$ then $N_\beta=(\Z_2^2\rtimes S_2)\times\Z_4$. Thus $N_\alpha/G=\Z_2^2$ and $N_\beta/G=\Z_2\times\Z_4$.  This example shows that even if $M_1$ and $M_2$ are diffeomorphic manifolds and they are quotients of the same complex torus $T$, their groups of automorphisms may be different.

Note that from the criterion given by \cite[Theorem 5.3]{S} (and Diagram \ref{diag:aff}), the group $\Aff(M)$ is infinite for any complex Hantzsche-Wendt manifold. The finiteness of the group $\Aut(M)$ is also a consequence of \cite[Theorem 0.1 (IV)]{OgSa} (see also \cite{LaOgPe}).

\subsection{Finiteness of \texorpdfstring{$\Aut M$}{Aut M}}
It is known that the finiteness of the group $\Aff(M)$ of a flat manifold $M$ depends only on its fundamental group (even more, on the $\Q$-isomorphism class of its holonomy representation, see \cite[Theorem 5.3]{S}). However for the automorphism group $\Aut(M)$ of a flat K\"ahler manifold situation is much different. 

Let us consider some examples which will be obtained by the following procedure. Let $M$ be a flat $n$-manifold with holonomy representation $\varrho$ defined by a cohomology class $\alpha\in H^2(G,\Z^n)$. For any subrepresentation $\sigma$ of $\varrho$ of degree $k$ there exists a flat $(n+k)$-manifold $M'$ with the same holonomy group and holonomy representation $\varrho\oplus\sigma$ (it may be defined by a class $(\alpha,0)\in H^2(G,\Z^n\oplus\Z^k)$). If $M$ is additionally K\"ahler, and $\sigma$ is an essentially complex subrepresentation, then $M'$ also is K\"ahler.

First we apply this construction to the CHW threefold $M$. This way we can get a flat K\"ahler manifold $M'$ with holonomy group $G=\langle a,b\rangle\simeq \Z_2^2$ with complex holonomy representation:
\[
a\mapsto \begin{bmatrix}1&&&\\&-1&&\\&&-1&\\&&&-1\end{bmatrix}, b\mapsto \begin{bmatrix}-1&&&\\&1&&\\&&-1&\\&&&-1\end{bmatrix}
\]
We may set $T=E_1\times E_2\times T_3$ in such a way that $E_1, E_2$ are elliptic curves, $T_3$ a complex $2$-torus and $\Aut_0(T)=\Aut_0(E_1)\times\Aut_0(E_2)\times\Aut_0(T_3)$. Then there exists $\tilde{G}\subset\Aut(T)$ such that underlying real manifold of $T/\tilde{G}$ is $M'$. Moreover $N_{\Aut_0(T)}(G)=\Aut_0(T)$, and in the consequence the finiteness of $\Aut(T/\tilde{G})$ depends only on the finiteness of $\Aut_0(T_3)$. As there exist complex $2$-tori with finite and infinite automorphism groups, we see that the finiteness of the group $\Aut(M')$ depends on the choice of the complex torus $T$.

Next example is obtained in similar way from the manifold with holonomy $\Z_3^2$. Consider a flat manifold with integral holonomy representation:
\[
\varrho(g_1)=
\begin{bmatrix}I&&&&\\&\Xi&&&\\&&\Xi&&\\&&&\Xi&\\&&&&\Xi\end{bmatrix},\ 
\varrho(g_2)=
\begin{bmatrix}\Xi&&&&\\&I&&&\\&&\Xi&&\\&&&\Xi^2&\\&&&&\Xi^2\end{bmatrix}
\]
(using the same notation as in Section \ref{sec:c3xc3}). Let $T=E_1\times E_2\times\ldots\times E_5$, where all curves $E_k$ are biholomorphic to $E_\xi$. The above holonomy representation may correspond to the following diagonal actions on $T$: $g_1=(1,\xi,\xi,\xi,\xi), g_2=(\xi,1,\xi,\xi^2,\xi^2)$ or $g_1'=(1,\xi,\xi,\xi,\xi^2), g_2'=(\xi,1,\xi,\xi^2,\xi)$. Groups $G=\langle g_1,g_2\rangle$ and $G'=\langle g_1',g_2'\rangle$ are not conjugated in $\Aut_0(T)$ (in fact, complex holonomy representations will be non-equivalent). Moreover the group $N_{\Aut_0(T)}(G)$ is infinite (as it contains any automorphism of $E_4\times E_5$) whereas $N_{\Aut_0(T)}(G')$ is finite (and isomorphic to $\Z_6^5\rtimes S$ where $S$ is a cyclic group generated by permutation $(1,2)(4,5)$). This time the finiteness of the automorphism group depends on the choice of the action of the holonomy group on $T$ (we may also say, on the complex holonomy representation).

\section{Open questions}
Let us end up with two questions, for which we hope to find an answer in the future.
\begin{enumerate}
	\item Does there exist a flat manifold $M$ with first Betti number equal to zero such that for any complex structure on $M$ the group $\Aut(M)$ is infinite?
	\item Is there a (necessarily incomplete) criterion for the finiteness of the group $\Aut(M)$ in terms of $\pi_1(M)$ or the complex holonomy representation of $M$?
\end{enumerate}

\section*{Acknowledgments}

The authors would like to thank Andrzej Szczepański for helpful discussions.

\end{document}